\documentclass[11pt]{article}

\usepackage{amssymb}
\usepackage{amsmath}
\usepackage{amsfonts}
\usepackage{epsfig}



\newcommand{\sect}[1]{\setcounter{equation}{0}\section{#1}}

\newcommand{\X}{\mathfrak{X}}
\newcommand{\Y}{\mathfrak{Y}}
\newcommand{\Conf}{\mathrm{Conf}}
\newcommand{\R}{\mathbb{R}}
\newcommand{\C}{\mathbb{C}}
\newcommand{\Z}{\mathbb{Z}}
\renewcommand{\P}{\mathcal{P}}
\newcommand{\T}{\mathcal{T}}
\newcommand{\virt}{\mathsf{virt}}
\newcommand{\ul}{\underline}
\newcommand{\LV}{\Lambda^{\frac\infty2}V}
\newcommand{\vac}{v_{\mathrm{vac}}}
\newcommand{\K}{\mathfrak{K}}

\begin{document}

\title{Determinantal point processes}

\author{Alexei Borodin
\footnote{Department of Mathematics, California Institute of
Technology, Pasadena, USA and IITP RAS, Moscow, Russia} }

\date{}

\maketitle

\begin{abstract}
 We present a list of algebraic, combinatorial, and analytic
mechanisms that give rise to determinantal point processes.
\end{abstract}

\setlength{\baselineskip}{5.0mm}

\newtheorem{theo}{Theorem}[section]

%\tableofcontents

\ \\

%%%%%%%%%%%%%%%%%%%%%%%%%%%%%%%%%%%%%%%%%%%%%%%%%%%%%%%%%%%%%%%%%%%%%%%%%%%
\sect{Introduction}\label{intro}

Let $\X$ be a discrete space. A (simple) random point process $\P$
on $\X$ is a probability measure on the set $2^\X$ of all subsets of
$\X$. $\P$ is called {\it determinantal\/} if there exists a
$|\X|\times|\X|$ matrix $K$ with rows and columns marked by elements
of $\X$, such that for any finite $Y=(y_1,\dots,y_n)\subset \X$ one
has
$$
\Pr\{X\in 2^\X\mid Y\subset X\}=\det[K(y_i,y_j)]_{i,j=1}^n.
$$
A similar definition can be given for $\X$ being any reasonable
space; then the measure lives on locally finite subsets of $\X$.

Determinantal point processes (with $\X=\R$) have been used in
random matrix theory since early 60's. As a separate class
determinantal processes were first singled out in \cite{Mac75} to
model fermions in thermal equilibrium, cf. \cite{Ben73}, and the
term `fermion' point processes was used.
 The term `determinantal' was introduced in
\cite{Bor00a}, for the reason that the particles of the process
studied there were of two kinds; particles of the same kind
repelled, while particles of different kinds attracted. Nowadays,
the expression `determinantal point process (or field)' is standard.

There are several excellent surveys of the subject available, see
\cite{Sos00}, \cite{Lyo03}, \cite{Joh05}, \cite{Kon05},
\cite{Hou06}, \cite{Sos06}. The reader may find there a detailed
discussion of probabilistic properties of determinantal processes as
well as a wide array of their applications; many applications are
also described in various chapters of this volume.

The goal of the present note is to bring together all known
algebraic, combinatorial, and analytic mechanisms that produce
determinantal processes. Many of the well-known determinantal
processes fit into more than one class described below. However,
none of the classes is superseded by any other.

\sect{Generalities}\label{general}

Let $\X$ be a locally compact separable topological space. A {\it
point configuration\/} $X$ in $\X$ is a locally finite collection of
points of the space $\X$. Any such point configuration is either
finite or infinite. For our purposes it suffices to assume that the
points of $X$ are always pairwise distinct. The set of all point
configurations in $\X$ will be denoted as $\Conf(\X)$.

A relatively compact Borel subset $A\subset \X$ is called {\it a
window}. For a window $A$ and $X\in\Conf(\X)$, set $N_A(X)=|A\cap
X|$ (number of points of $X$ in the window). Thus, $N_A$ can be
viewed as a function on $\Conf(\X)$. We equip $\Conf(\X)$ with the
Borel structure generated by functions $N_A$ for all windows $A$.

A {\it random point process\/} on $\X$ is a probability measure on
$\Conf(\X)$.

Given a random point process, one can usually define a sequence
$\{\rho_n\}_{n=1}^\infty$, where $\rho_n$ is a symmetric measure on
$\X^n$ called the $n$th {\it correlation measure}. Under mild
conditions on the point process, the correlation measures exist and
determine the process uniquely, cf. \cite{Len73}.

The correlation measures are characterized by the following
property: For any $n\ge1$ and a compactly supported bounded Borel
function $f$ on $\X^n$ one has
\begin{equation}\label{corr1}
\int_{\X^n}f\rho_n=\left\langle \sum_{x_{i_1},\dots,x_{i_n}\in
X}f(x_{i_1},\dots,x_{i_n})\right \rangle_{X\in\Conf(\X)}
\end{equation}
where the sum on the right is taken over all $n$-tuples of pairwise
distinct points of the random point configuration $X$.

Often one has a natural measure $\mu$ on $\X$ (called the {\it
reference measure\/}) such that the correlation measures have
densities with respect to $\mu^{\otimes n}$, $n=1,2,...\,$. Then the
density of $\rho_n$ is called the $n$th {\it correlation function\/}
and it is usually denoted by the same symbol ``$\rho_n$''.

If $\X\subset \R$ and $\mu$ is absolutely continuous with respect to
the Lebesgue measure, then the probabilistic meaning of the $n$th
correlation function is that of the density of probability to find
an eigenvalue in each of the infinitesimal intervals around points $
x_1, x_2, \ldots x_n $:
\begin{multline*}
\rho_n(x_1, x_2, \ldots
x_n) \mu(dx_1) \cdots \mu(dx_n)\\= \Pr \,\{ \text{there is a
particle in each interval}\ \ (x_i, x_i +dx_i)\}.
\end{multline*}
On the other hand, if $\mu$ is supported by a discrete set of
points, then
\begin{multline*}
 \rho_n(x_1, x_2, \ldots x_n) \mu(x_1) \cdots \mu(x_n)\\ = \Pr
\{ \text{there is a particle at each of the points}\ \ x_i\}.
\end{multline*}

Assume that we are given a point process $\P$ and a reference
measure such that all correlation functions exist. The process $\P$
is called {\it determinantal\/} if there exists a function
$K:\X\times\X\to \C$ such that
\begin{equation}\label{corr}
\rho_n(x_1,\dots,x_n)=\det[K(x_i,x_j)]_{i,j=1}^n,\qquad
n=1,2,\dots\,.
\end{equation}
The function $K$ is called a {\it correlation kernel} of $\P$.

The determinantal form of the correlation functions (\ref{corr})
implies that many natural observables for $\P$ can be expressed via
the kernel $K$. We mention a few of them. For the sake of
simplicity, we assume that the state space $\X$ is discrete and
$\mu$ is the counting measure; under appropriate assumptions, the
statements are easily carried over to more general state spaces.

\begin{itemize}

\item Let $I$ be a (possibly infinite) subset of $\X$. Denote by
$K_I$ the operator in $\ell^2(I)$ obtained by restricting the kernel
$K$ to $I$. Assume that $K_I$ is a trace class operator.
\footnote{For discrete $\X$, a convenient sufficient condition for
$K_I$ to be of trace class is $ \sum_{x,y\in I} |K(x,y)|<\infty$.}
Then the intersection of the random configuration $X$ with $I$ is
finite almost surely and
\begin{equation*}
\Pr\{|X\cap I|=N\}=\frac{(-1)^N}{N!}\left.\frac{d^N}{dz^N}\det\Big(
\mathbf{1}-zK_I\Big)\right|_{z=1}\,.
\end{equation*}
In particular, the probability that $X\cap I$ is empty is equal to
\begin{equation*}
\Pr\{X\cap I=\varnothing\}=\det\Big(\mathbf{1}-K_I\Big)\,.
\end{equation*}

More generally, if $I_1,\dots,I_m$ is a finite family of pairwise
nonintersecting intervals such that the operators
$K_{I_1},\dots,K_{I_m}$ are trace class then
\begin{multline}\label{e41}
\Pr\{|X\cap I_1|=N_1,\dots,|X\cap I_m|=N_m \}\\
=\left.\frac{(-1)^{\sum_{i=1}^m N_i}}{\prod_{i=1}^m
N_i!}\frac{\partial^{N_1+\dots+N_m}} {\partial
z_1^{N_1}\dots\partial z_m^{N_m}}
\det\Big(\mathbf{1}-z_1K_{I_1}-\dots-z_mK_{I_m}\Big)\right|_{z_1=\dots=z_m=1}.
\end{multline}
\item Slightly more generally, let $\phi$ be a function on $\X$ such that
the kernel $(1-\phi(x))K(x,y)$ defines a trace class operator
$(1-\phi)K$ in $\ell^2(\X)$. Then
\begin{equation}\label{gen}
\mathbb{E}\, \left(\prod_{x_i\in X} \phi(x_i)\right)
=\det(\mathbf{1}-(1-\phi)K).
\end{equation}
Specifying $\phi=\sum_{j=1}^m (1-z_j) \mathbf{1}_{I_j}$ leads to
(\ref{e41}).

\item For $I\subset \X$ such that $K_I$ is trace class and $\det(\mathbf{1}-K_I)\ne 0$,
 and arbitrary pairwise distinct locations
$\{x_1,\dots,x_n\}\subset I$, $n=1,2,\dots$, set
\begin{multline*}
\mathcal{J}_{I,n}(x_1,\dots,x_n)=\Pr\{ \text{there is a particle at
each of the points $x_i$}\\  \text{and there are no other particles
in $I$}\}.
\end{multline*}
These are sometimes called {\it Janossy measures\/}. One has
\begin{equation}\label{janossy}
\mathcal{J}_{I,n}(x_1,\dots,x_n)=\det(\mathbf{1}-K_I)\cdot
\det[L_I(x_i,x_j)]_{i,j=1}^n,
\end{equation}
where $L_I$ is the matrix of the operator
$K_I(\mathbf{1}-K_I)^{-1}$.
\end{itemize}

Simple linear-algebraic proofs of (\ref{gen}) and (\ref{janossy})
can extracted from the proof of Proposition A.6 in \cite{Bor00b}. We
also refer to Chapter 4 in this volume for a detailed discussion of
(\ref{e41})--(\ref{janossy}) and many related identities.

\sect{Loop-free Markov chains}

Let $\X$ be a discrete space, and let $P=[P_{xy}]_{x,y\in\X}$ be the
matrix of transition probabilities for a discrete time Markov chain
on $\X$. That is, $P_{xy}\ge 0$ for all $x,y\in\X$ and
$$
\sum_{y\in\X} P_{xy}=1 \quad \text{for any}\quad  x\in\X.
$$

Let us assume that our Markov chain is {\it loop-free}, i.e. the
trajectories of the Markov chain do not pass through the same point
twice almost surely. In other words, we assume that
$$
(P^k)_{xx}=0\qquad \text{for any}\quad k>0 \quad\text{and}\quad
x\in\X.
$$

This condition guarantees the finiteness of the matrix elements of
the matrix
$$
Q=P+P^2+P^3+\dots\,.
$$
Indeed, $(P^k)_{xy}$ is the probability that the trajectory started
at $x$ is at $y$ after $k$th step. Hence, $Q_{xy}$ is the
probability that the trajectory started at $x$ passes through $y\ne
x$, and since there are no loops we have $Q_{xy}\le 1$. Clearly,
$Q_{xx}\equiv 0$.

The following (simple) fact was proved in \cite{Bor08a}.

\begin{theo}\label{markov} For any probability measure
$\pi=[\pi_x]_{x\in\X}$ on $\X$, consider the Markov chain with
initial distribution $\pi$ and transition matrix $P$ as a
probability measure on trajectories viewed as subsets of $\X$. Then
this measure on $2^\X$ is a determinantal point process on $\X$ with
correlation kernel
$$
K(x,y)=\pi_x+(\pi Q)_x-Q_{yx}.
$$
\end{theo}

Note that the correlation kernel is usually not
self-adjoint\footnote{In fact, it can be written as a sum of a
nilpotent matrix and a matrix of rank 1.}, and self-adjoint examples
should be viewed as ``exotic''. One such example goes back to
\cite{Mac75}, see also \S2.4 of \cite{Sos00}: It is a 2-parameter
family of renewal processes
--- processes on $\Z$ or $\R$ with positive i.i.d. increments.
Theorem \ref{markov} implies that if we do not insist on
self-adjointness then any process with positive i.i.d. increments is
determinantal.

\sect{Measures given by products of determinants}

Let $\X$ be a finite set and $N$ be any natural number no greater
than $|\X|$. Let $\Phi_n$ and $\Psi_n$, $n=1,2,\dots,N$, be
arbitrary complex-valued functions on $\X$. To any point
configuration $X\in \Conf(\X)$ we assign its weight $W(X)$ as
follows: If the number of points in $X$ is not $N$ then $W(X)=0$.
Otherwise, using the notation $X=\{x_1,\dots,x_N\}$, we have
$$
W(X)=\det\left[\Phi_i(x_j)\right]_{i,j=1}^N\det\left[\Psi_i(x_j)\right]_{i,j=1}^N.
$$

Assume that the partition function of our weights does not vanish
$$
Z:=\sum_{X\in\Conf(\X)} W(X)\ne 0.
$$
Then the normalized weights $\widetilde W(X)=W(X)/Z$ define a
(generally speaking, complex valued) measure on $\Conf(\X)$ of total
mass 1. Such measures are called {\it biorthogonal
ensembles\/}.\footnote{This term was introduced in \cite{Bor99} and
is now widely used.} For complex valued point processes we use
(\ref{corr1}) to define their correlation functions.

An especially important subclass of biorthogonal ensembles consists
of \emph{orthogonal polynomial ensembles}, for which $\X$ must be a
subset of $\C$, and
$$
W(X)=\prod_{1\le i<j\le N} |x_i-x_j|^2\cdot\prod_{i=1}^N w(x_i)
$$
for a function $w:\X\to\R_+$, see e.g. \cite{Kon05} and Chapter 4 of
this volume.

\begin{theo}\label{bi} Any biorthogonal ensemble is a determinantal point
process. Its correlation kernel has the form
$$
K(x,y)=\sum_{i,j=1}^N \left[G^{-t}\right]_{ij}\Phi_i(x)\Psi_j(y),
$$
where $G=[G_{ij}]_{i,j=1}^N$ is the Gram matrix:
$G_{ij}=\sum_{x\in\X}\Phi_i(x)\Psi_j(x)$.\footnote{The invertibility
of the Gram matrix is implied by the assumption $Z\ne 0$.}
\end{theo}

The statement immediately carries over to $\X$ being an arbitrary
state space with reference measure $\mu$; then one has
$G_{ij}=\int_{\X}\Phi_i(x)\Psi_j(x)\mu(dx)$.

Probably the first appearance of Theorem \ref{bi} is in the seminal
work of F.~J.~Dyson \cite{Dys62a}, where it was used to evaluate the
correlation functions of the eigenvalues of the Haar-distributed
$N\times N$ unitary matrix. In that case, $\X$ is the unit circle,
$\mu$ is the Lebesgue measure on it,
$$
\Phi_i(z)=z^{i-1},\qquad \Psi_i(z)=\bar{z}^{i-1},\qquad |z|=1,\quad
i=1,\dots,N,
$$
and the Gram matrix $G$ coincides with the identity matrix.

In the same volume, Dyson \cite{Dys62b} introduced a Brownian motion
model for the eigenvalues of random matrices (currently known as the
{\it Dyson Brownian motion}), and it took more than three decades to
find a determinantal formula for the time-dependent correlations of
eigenvalues in the unitarily invariant case. The corresponding claim
has a variety of applications; let us state it. Again, for
simplicity of notation, we work with finite state spaces.

Let $\X^{(1)},\dots,\X^{(k)}$ be finite sets. Set
$\X=\X^{(1)}\sqcup\dots\sqcup\X^{(k)}$. Fix a natural number $N$.
Let
\begin{gather*}
\Phi_i:\X^{(1)}\to\C,\qquad  \Psi_i:\X^{(k)}\to\C, \qquad
i=1,\dots,N \\
\T_{j,j+1}:\X^{(j)}\times \X^{(j+1)}\to \C,\qquad j=1,\dots,k-1,
\end{gather*}
be arbitrary functions. To any $X\in\Conf(\X)$ assign its weight
$W(X)$ as follows. If $X$ has exactly $N$ points in each $X^{(j)}$,
$j=1,\dots,k$ then denoting $X\cap \X^{(j)}=\{x_1^{(j)},\dots,
x_N^{(j)}\}$ we have
\begin{multline}\label{EMweight}
W(X)=\det\left[\Phi_i(x_j^{(1)})\right]_{i,j=1}^N
\det\left[\T_{1,2}(x_i^{(1)},x_j^{(2)})\right]_{i,j=1}^N\cdots\\
\times \det\left[\T_{k-1,k}(x_i^{(k-1)},x_j^{(k)})\right]_{i,j=1}^N
\det\left[\Psi_i(x_j^{(1)})\right]_{i,j=1}^N;
\end{multline}
otherwise $W(X)=0$.

As for biorthogonal ensembles above, we assume that the partition
function of these weights is nonzero and define the corresponding
normalized set of weights. This gives a (generally speaking, complex
valued) random point process on $\X$.

In what follows we use the notation
\begin{gather*}
(f*g)(x,y)=\sum_z f(x,z)g(z,y),\quad h_1*h_2=\sum_x h_1(x)h_2(x),\\
(h_1*f)(y)=\sum_x h_1(x)f(x,y),\quad (g*h_2)(x)=\sum_y g(x,y)h_2(y)
\end{gather*}
for arbitrary functions $f(x,y)$, $g(x,y)$, $h_1(x)$, $h_2(x)$,
where the sums are taken over all possible values of the summation
variables.

\begin{theo}\label{EM}
The random point process defined by (\ref{EMweight}) is
determinantal. The correlation kernel on $\X^{(p)}\times\X^{(q)}$,
$p,q=1,\dots,N$, can be written in the form
\begin{multline}\label{EMkernel}
K(x^{(p)},y^{(q)})=-\mathbf{1}_{p>q}\cdot(\T_{q,q+1}*\dots*\T_{p-1,p})(y^{(q)},x^{(p)})\\+
\sum_{i,j=1}^N \left[G^{-t}\right]_{ij}
\left(\Phi_i*\T_{1,2}*\dots*\T_{p-1,p}\right)(x^{(p)})
\left(\T_{q,q+1}*\dots*\T_{k-1,k}*\Psi_j\right)(y^{(q)}),
\end{multline}
where the Gram matrix $G=\left[G_{ij}\right]_{i,j=1}^N$ is defined
by
$$
G_{ij}=\Phi_i*\T_{1,2}*\dots*\T_{k-1,k}*\Psi_j,\qquad i,j=1,\dots,N.
$$
\end{theo}

Similarly to Theorem \ref{bi}, the statement is easily carried over
to general state spaces $\X^{(j)}$.

Theorem \ref{EM} is often referred to as the {\it Eynard-Mehta
theorem}, it was proved in \cite{Eyn98} and also independently in
\cite{Nag98}. Other proofs can be found in \cite{Joh03},
\cite{Tra04}, \cite{Bor05}.

The algebraically ``nice'' case of the Eynard-Mehta theorem, which
e.g. takes place for the Dyson Brownian motion, consists in the
existence of an orthonormal basis $\{\Xi^{(j)}_i\}_{i\ge 1}$ in each
$L^2(\X^{(j)})$, $j=1,\dots,k$, such that
$$
T_{j,j+1}(x,y)=\sum_{i\ge 1}c_{j,j+1;i}\,\Xi^{(j)}_i(x)
\Xi^{(j+1)}_i(y),\qquad j=1,2,\dots,k-1,
$$
for some constants $c_{j,j+1;i}$, and
\begin{gather*}
\mathrm{Span}\{\Xi^{(1)}_1,\dots,\Xi^{(1)}_N\}=\mathrm{Span}\{\Phi_1,\dots,\Phi_N\},\\
\mathrm{Span}\{\Xi^{(k)}_1,\dots,\Xi^{(k)}_N\}=\mathrm{Span}\{\Psi_1,\dots,\Psi_N\}.
\end{gather*}
Then, with the notation $c_{k,l;i}=c_{k,k+1;i}c_{k+1,k+2;i}\cdots
c_{l-1,l;i}$, (\ref{EMkernel}) reads
$$
K(x^{(p)},y^{(q)})=\begin{cases} \displaystyle\sum_{i=1}^N
\frac{1}{c_{p,q;i}}\,\Xi^{(p)}_i(x^{(p)})\,\Xi^{(q)}_i(y^{(q)}),&p\le
q,\\
\displaystyle  -\sum_{i>N}
{c_{q,p;i}}\,\Xi^{(p)}_i(x^{(p)})\,\Xi^{(q)}_i(y^{(q)}),&p>q.
\end{cases}
$$

The ubiquitousness of the Eynard-Mehta theorem in applications is
explained by the combinatorial statement known as the
Lindstr\"om-Gessel-Viennot (LGV) theorem, see \cite{Ste90} and
references therein, that we now describe.

Consider a finite\footnote{The assumption of finiteness is not
necessary as long as the sums in (\ref{paths}) converge.} directed
acyclic graph and denote by $V$ and $E$ the sets of its vertices and
edges. Let $w:E\to\C$ be an arbitrary weight function. For any path
$\pi$ denote by $w(\pi)$ the product of weights over the edges in
the path: $w(\pi)=\prod_{e\in \pi} w(e)$. Define the weight of a
collection of paths as the product of weights of the paths in the
collection (we will use the same letter $w$ to denote it). We say
that two paths $\pi_1$ and $\pi_2$ do not intersect (notation
$\pi_1\cap\pi_2=\varnothing$) if they have no common vertices.

For any $u,v\in V$, let $\Pi(u,v)$ be the set of all (directed)
paths from $u$ to $v$. Set
\begin{equation}\label{paths}
\mathcal{T}(u,v)=\sum_{\pi\in\Pi(u,v)} w(\pi).
\end{equation}

\begin{theo}\label{LGV} Let $(u_1,\dots,u_n)$ and $(v_1,\dots,v_n)$ be two
$n$-tuples of vertices of our graph, and assume that for any
nonidentical permutation $\sigma\in S(n)$,
$$
\left\{(\pi_1,\dots,\pi_n)\mid \pi_i\in
\Pi\left(u_i,v_{\sigma(i)}\right),\ \pi_i\cap\pi_j=\varnothing,\
i,j=1,\dots,n\right\}=\varnothing.
$$
Then
$$
\sum_{\substack{\pi_1\in \Pi(u_1,v_1),\dots, \pi_n\in \Pi(u_n,v_n)\\
\pi_i\cap \pi_j=\varnothing,\
i,j=1,\dots,n}}w(\pi_1,\dots,\pi_n)=\det\left[\mathcal{T}(u_i,v_j)\right]_{i,j=1}^n.
$$
\end{theo}

Theorem \ref{LGV} means that if, in a suitable weighted oriented
graph, we have nonintersecting paths with fixed starting and ending
vertices, then the distributions of the intersection points of these
paths with any chosen ``sections'' have the same structure as
(\ref{EMweight}), and thus by Theorem \ref{EM} we obtain a
determinantal point process.

A continuous time analog of Theorem \ref{LGV} goes back to
\cite{Kar59}, who in particular proved the following statement (the
next paragraph is essentially a quotation).

Consider a stationary stochastic process whose state space is an
interval on the extended real line. Assume that the process has
strong Markov property and that its paths are continuous everywhere.
Take $n$ points $x_1<\dots<x_n$ and $n$ Borel sets $E_1<\dots<E_n$,
and suppose $n$ labeled particles start at $x_1,\dots,x_n$ and
execute the process simultaneously and independently. Then the
determinant $\det\left[P_t(x_i,E_j)\right]_{i,j=1}^n$, with
$P_t(x,E)$ being the transition probability of the process, is equal
to the probability that at time $t$ the particles will be found in
sets $E_1,\dots,E_n$ respectively without any of them ever having
been coincident in the intervening time.

Similarly to Theorem \ref{LGV}, this statement coupled with Theorem
\ref{EM} leads to determinantal processes, and this is exactly the
approach that allows one to compute the time-dependent eigenvalue
correlations of the Dyson Brownian motion.

We conclude this section with a generalization of the Eynard-Mehta
theorem that allows the number of particles to vary.

Let $\X_1,\dots,\X_N$ be finite sets, and
\begin{equation*}
\begin{aligned}
\phi_n(\,\cdot\,,\,\cdot\,):\X_{n-1}\times\X_{n}\to \C,\qquad  &n=2,\dots,N, \\
\phi_n(\virt,\,\cdot\,):\X_{n}\to\C,\qquad  &n=1,\dots,N,\\
\Psi_j(\,\cdot\,):\X_N\to\C,\qquad &j=1,\dots,N,
\end{aligned}
\end{equation*}
be arbitrary functions on the corresponding sets. Here the symbol
$\virt$ stands for a ``virtual'' variable, which is convenient to
introduce for notational purposes. In applications, $\virt$ can
sometimes be replaced by $+\infty$ or $-\infty$.

Let $c(1),\dots,c(N)$ be arbitrary nonnegative integers, and let
\begin{equation*}
t_{0}^N\leq\dots\leq t_{c(N)}^N= t_{0}^{N-1}\leq\dots\leq
t_{c(N-1)}^{N-1}=t_0^{N-2}\leq \dots\leq t^2_{c(2)}=t^1_0\leq\dots
\leq t^1_{c(1)}
\end{equation*}
be real numbers. In applications, these numbers may refer to time
moments of an associated Markov process. Finally, let
\begin{equation*}
{\cal T}_{t_a^n,t_{a-1}^n}(\,\cdot\,,\,\cdot\,):\X_n\times\X_n\to
\C,\qquad n=1,\dots,N,\quad a=1,\dots,c(n),
\end{equation*}
be arbitrary functions.

Set $\X=(\X_1\sqcup\dots\sqcup\X_1)\sqcup\dots\sqcup(\X_N\sqcup
\dots\sqcup\X_N)$ with $c(n)+1$ copies of each
$\X_n$\footnote{Instead of $c(n)+1$ copies of $\X_n$ one can take
same number of different spaces, and a similar result will hold. We
decided not to do it in order not to clutter the notation anymore.},
and to any $X\in\Conf(\X)$ assign its weight $W(X)$ as follows.

The weight $W(X)$ is zero unless $X$ has exactly $n$ points in each
copy of $\X_n$, $n=1,\dots,N$. In the latter case, denote the points
of $X$ in the $m$th copy of $\X_n$ by $x^n_k(t^n_m)$, $k=1,\dots,n$,
and set
\begin{equation}\label{GMgen}
\begin{aligned}
W(X)= & \prod_{n=1}^{N} \Bigg[\det{\Bigl[\phi_n\bigl(x_k^{n-1}(t_0^{n-1}),x_l^n(t^n_{c(n)})\bigr)\Bigr]}_{k,l=1}^n \\
&\ \  \times\prod_{a=1}^{c(n)} \det{\Bigl[{\cal
T}_{t_a^n,t_{a-1}^n}\bigl(x_k^n(t^n_a),x^n_l(t^n_{a-1})\bigr)\Bigr]}_{k,l=1}^n
\Bigg]\cdot\det{\Bigl[\Psi_{l}\bigl(x^{N}_k(t_0^{N})\bigr)\Bigr]}_{k,l=1}^N,
\end{aligned}
\end{equation}
where $x^{n-1}_{n}(\,\cdot\,)=\virt$ for all $n=1,\dots,N$.

Once again, we assume that the partition function does not vanish,
and normalizing the weights we obtain a (generally speaking, complex
valued) point process on $\X$.

We need more notation. For any $n=1,\dots,N$ and two time moments
$t_a^n>t_b^n$ we define
\begin{equation*}
{\cal T}_{t_a^n,t_b^n}={\cal T}_{t_a^n,t_{a-1}^n}*{\cal
T}_{t_{a-1}^n,t_{a-2}^n}*\cdots *{\cal T}_{t_{b+1}^n,t_b^n},\qquad
{{\cal T}}^n = {{\cal T}}_{t_{c(n)}^n,t_0^n}.
\end{equation*}
 For any
time moments $t_{a_1}^{n_1}\ge t_{a_2}^{n_2}$ with $(a_1,n_1)\ne
(a_2,n_2)$, we denote the convolution over all the transitions
between them by $\phi^{(t_{a_1}^{n_1},t_{a_2}^{n_2})}$:
\begin{equation*}
\phi^{(t_{a_1}^{n_1},t_{a_2}^{n_2})}={{\cal
T}}_{t^{n_1}_{a_1},t^{n_1}_{0}} * \, \phi_{n_1+1}*{{\cal T}}^{n_1+1}
*\cdots*\phi_{n_2}*{{\cal T}}_{t^{n_2}_{c(n_2)},t^{n_2}_{a_2}}.
\end{equation*}
If there are no such transitions, i.~e. if
$t_{a_1}^{n_1}<t_{a_2}^{n_2}$ or $(a_1,n_1)=(a_2,n_2)$, we set
$\phi^{(t_{a_1}^{n_1},t_{a_2}^{n_2})}=0$.

Furthermore, define the ``Gram matrix'' $G={\left[
G_{kl}\right]}_{k,l=1}^N$ by
\begin{equation*}
G_{kl}=\big(\phi_{k}*{{\cal T}}^k*\cdots *\phi_{N}*{{\cal
T}}^{N}*\Psi_{l}\big)(\virt),\qquad k,l=1,\dots,N,
\end{equation*}
and set
\begin{equation*}
\Psi^{t^n_a}_{l}=\phi^{(t^n_a,t^{N}_0)}*\Psi_{l},\qquad l=1,\dots,N.
\end{equation*}

\begin{theo}\label{ThmPushASEP}
The random point process on $\X$ defined by (\ref{GMgen}) is
determinantal. Its correlation kernel can be written in the form
\begin{equation*}
\begin{aligned}
K(t^{n_1}_{a_1},x_1; t^{n_2}_{a_2},x_2)&= -\phi^{(t^{n_2}_{a_2},t^{n_1}_{a_1})}(x_2,x_1) \\
+&  \sum_{i=1}^{n_1}\sum_{j=1}^{N}\left[G^{-t}\right]_{ij}(\phi_i *
\phi^{(t^i_{c(i)},t^{n_1}_{a_1})})(\virt,x_1)\,
\Psi^{t^{n_2}_{a_2}}_{j}(x_2)
 .
\end{aligned}
\end{equation*}
\end{theo}

One proof of Theorem~\ref{ThmPushASEP} was given in~\cite{Bor08b};
another proof can be found in Section 4.4 of~\cite{For08}. Although
we stated Theorem~\ref{ThmPushASEP} for the case when all sets
$\X_n$ are finite, one easily extends it to a more general setting.

\sect{L-ensembles}

The definition of L-ensembles is closely related to (\ref{janossy}).

Let $\X$ be a finite set.
 Let $L$ be a $|\X|\times|\X|$
matrix whose rows and column are parameterized by points of $\X$.
For any subset $X\subset \X$ we will denote by $L_X$ the symmetric
submatrix of $L$ corresponding to $X$: $
L_X=\left[L(x_i,x_j)\right]_{x_i,x_j\in X}.$ If determinants of all
such submatrices are nonnegative (e.g., if $L$ is positive
definite), one can define a random point process on $\X$ by
$$
\Pr\{X\}=\frac{\det L_X}{\det(\mathbf{1}+L)}\,,\qquad X\subset \X.
$$
This process is called the {\it $L$-ensemble}.

The following statement goes back to \cite{Mac75}.

\begin{theo}\label{L} The $L$-ensemble as defined
above is a determinantal point process with the correlation kernel
$K$ given by $K=L(\mathbf{1}+L)^{-1}$.
\end{theo}

Take a nonempty subset $\Y$ of $\X$ and, given an $L$-ensemble on
$\X$, define a new random point process on $\Y$ by considering the
intersections of the random point configurations $X\subset \X$ of
the $L$-ensemble with $\Y$, provided that these point configurations
contain the complement $\overline \Y$ of $\Y$ in $\X$. It is not
hard to see that this new process can be defined by
$$
\Pr \{Y\}=\frac{\det
L_{Y\cup\overline\Y}}{\det(\mathbf{1}_\Y+L)}\,,\qquad Y\in\Conf(\Y).
$$
Here $\mathbf{1}_\Y$ is the block matrix $\bmatrix
\mathbf{1}&0\\0&0\endbmatrix$ where the blocks correspond to the
splitting $\X=\Y\sqcup\overline{\Y}$. This new process is called the
{\it conditional $L$-ensemble}. The next statement was proved in
\cite{Bor05}.

\begin{theo}\label{condL} The conditional $L$-ensemble is a determinantal point
process with the correlation kernel given by
$$
K=\mathbf{1}_\Y-(\mathbf{1}_\Y+L)^{-1}\bigl|_{\Y\times\Y}.
$$
\end{theo}

Note that for $\Y=\X$, Theorem \ref{condL} coincides with Theorem
\ref{L}.

Not every determinantal process is an L-ensemble; for example, the
processes afforded by Theorem \ref{bi} have exactly $N$ particles,
which is not possible for an L-ensemble. However, as shown in
\cite{Bor05}, every determinantal process (on a finite set) is a
conditional L-ensemble.

The definition of L-ensembles and Theorems \ref{L}, \ref{condL} can
be carried over to infinite state spaces $\X$, given that $L$
satisfies appropriate conditions. In particular, the Fredholm
determinant $\det(\mathbf{1}_\Y+L)$ needs to be well defined.

Although L-ensembles do arise naturally, see e.~g.~\cite{Bor00b},
they also constitute a convenient computation tool. For example,
proofs of Theorems \ref{EM} and \ref{ThmPushASEP} given in
\cite{Bor05} and \cite{Bor08b} represent the processes in question
as conditional L-ensembles and employ Theorem \ref{condL}.

Here is another application of Theorem \ref{condL}.

A random point process on (a segment of) $\Z$ is called
\emph{one-dependent} if for any two finite sets $A,B\subset \Z$ with
$dist(A,B)\ge 2$, the correlation function factorizes:
$\rho_{|A|+|B|}(A\cup B)=\rho_{|A|}(A)\rho_{|B|}(B).$

\begin{theo} \label{rdet} Any one-dependent point process on (a segment
of) $\mathbb{Z}$ is determinantal. Its correlation kernel can be
written in the form
\[  K(x,y)=\begin{cases}{}
0, & x-y \geq 2,\\
-1, & x-y = 1,\\
\sum\limits_{r=1}^{y-x+1} (-1)^{r-1}
\sum\limits_{x=l_0<l_1<\cdots<l_r=y+1}  R_{l_0,l_1}R_{l_1,l_2}
\cdots R_{l_{r-1},l_r}& x\leq y,
\end{cases} \]
where $R_{a,b}=\rho_{b-a}(a,a+1,\dots,b-1)$.
\end{theo}

Details and applications can be found in \cite{Bor09}.

\sect{Fock space}

A general construction of determinantal point processes via the Fock
space formalism can be quite technical, see e.~g.~\cite{Lyt02}, so
we will consider a much simpler (however nontrivial) example
instead.

Recall that a \emph{partition}
$\lambda=(\lambda_1\ge\lambda_2\ge\dots\ge 0)$ is a weakly
decreasing sequence of nonnegative integers with finitely many
nonzero terms. We will use standard notations
$|\lambda|=\lambda_1+\lambda_2+\dots$ for the size of the partition
and $\ell(\lambda)$ for the number of its nonzero parts.

The \emph{poissonized Plancherel measure} on partitions is defined
by
\begin{equation}\label{pp}
\Pr\{\lambda\}=e^{-\theta^2}\left(\frac{\prod_{1\le i<j\le L}
(\lambda_i-i-\lambda_j+j)}{\prod_{i=1}^L(\lambda_i-i+L)!}\,\theta^{|\lambda|}\right)^2,
\end{equation}
where $\theta>0$ is a parameter, and $L$ is an arbitrary integer
$\ge \ell(\lambda)$. We refer to \cite{Bor00b}, \cite{Joh01} and
references therein for details.

It is convenient to parameterize partitions by subsets of
$\Z':=\Z+\frac 12$:
$$
\lambda\mapsto \mathcal{L}(\lambda)=\left\{\lambda_i-i+\tfrac
12\right\}_{i\ge 1}\subset \Z'.
$$
The pushforward of (\ref{pp}) via $\mathcal{L}$ defines a point
process on $\Z'$, and we aim to show that it is determinantal. We
follow \cite{Oko01}; other proofs of this fact can be found in
\cite{Bor00b}, \cite{Joh01}.

Let $V$ be a linear space with basis $\left\{\ul{k}\mid
k\in\Z'\right\}$. The linear space $\LV$  is, by definition, spanned
by vectors
$$
v_S=\ul{s_1} \wedge \ul{s_2} \wedge  \ul{s_3} \wedge  \dots\,,
$$
where $S=\{s_1>s_2>\dots\}\subset \Z'$ is such that both sets $S_+ =
S \setminus \Z'_{<0}$ and $S_- = \Z'_{<0} \setminus S $ are finite.
We equip $\LV$ with the inner product in which the basis $\{v_S\}$
is orthonormal.

\emph{Creation} and \emph{annihilation} operators in $\LV$ are
introduced as follows. The creation operator $\psi_k$ is the
exterior multiplication by $\ul{k}$: $\psi_k \left(f\right) = \ul{k}
\wedge f \,.$ The annihilation operator $\psi^*_k$ is its adjoint.
These operators satisfy the canonical anti-commutation relations
$$
\psi_k \psi^*_l + \psi^*_l \psi_k = \delta_{k,l}\,, \qquad
k,l\in\Z'.
$$
Observe that
\begin{equation}\label{psi}
\psi_k \psi^*_k \,\, v_S =
\begin{cases}
v_S\,, & k \in S \,, \\
0 \,, & k \notin S \,.
\end{cases}
\end{equation}

Let $C$ be the \emph{charge operator}: $Cv_S=(|S_+|-|S_-|)v_S$. One
easily sees that the zero-charge subspace $\ker C\subset \LV$ is
spanned by the vectors $v_{\mathcal{L}(\lambda)}$ with $\lambda$
varying over all partitions. The \emph{vacuum vector}
$$
\vac=\ul{-\tfrac12} \wedge \ul{-\tfrac32} \wedge \ul{-\tfrac52}
\wedge \dots
$$
corresponds to the partition with no nonzero parts.

Define the operators $\alpha_n=\sum_{k\in\Z'}\psi_{k-n}\psi^*_k$,
$n\in\Z\setminus\{0\}$. Although the sums are infinite, the
application of $\alpha_n$ to any $v_S$ yields a finite linear
combination of basis vectors. These operators satisfy the Heisenberg
commutation relations
$$
\alpha_m \alpha_n-\alpha_n\alpha_m = m \, \delta_{n,-m} \,,\qquad
m,n\in\Z\setminus\{0\}.
$$

For any $\theta>0$, define
$\Gamma_{\pm}(\theta)=\exp(\theta\alpha_{\pm 1})$. It is not
difficult to show that
\begin{equation}\label{gamma}
\Gamma_{\pm}^*(\theta)=\Gamma_{\mp}(\theta),\qquad
\Gamma_+(\theta)\Gamma_-(\theta')=e^{\theta\theta'}\cdot\Gamma_-(\theta')\Gamma_+(\theta),
\qquad \Gamma_+(\theta)\vac=\vac.
\end{equation}
One also proves that
$$
\Gamma_-(\theta)\vac=\sum_\lambda \left(\frac{\prod_{1\le i<j\le L}
(\lambda_i-i-\lambda_j+j)}{\prod_{i=1}^L(\lambda_i-i+L)!}\,\theta^{|\lambda|}\right)
v_{\mathcal{L}(\lambda)},
$$
where the sum is taken over all partitions, cf. (\ref{pp}). This
implies, together with (\ref{psi}), that for any $n\ge 1$ and
$x_1,\dots,x_n\in\Z'$, the correlation function of our point process
can be written as a matrix element
$$
\rho_n(x_1,\dots,x_n)=e^{-\theta^2}\left( \left( \prod_{i=1}^n
\psi_{x_i} \psi^*_{x_i} \right) \Gamma_-(\theta) \,
\vac,\Gamma_-(\theta)\vac\right).
$$
Using (\ref{gamma}) we obtain
\begin{equation}\label{e11}
\rho_n(x_1,\dots,x_n)=\left( \prod_{i=1}^n \Psi_{x_i} \Psi^*_{x_i}
\, \vac,\vac\right)\,,
\end{equation}
where
$$
\Psi_k = G\, \psi_k \, G^{-1} \,, \quad \Psi^*_k = G\, \psi^*_k \,
G^{-1} \,, \quad G= \Gamma_+(\theta)\, \Gamma_-(\theta)^{-1}\,.
$$
\begin{theo}\label{wedge} We have
\begin{equation}\label{e12}
\rho_n(x_1,\dots,x_n)=\det\left[K(x_i,x_j)\right]_{i,j=1}^n \,,
\end{equation}
where $K(x,y)= \left( \Psi_x \Psi^*_y \, \vac,\vac\right)$.
\end{theo}

The passage from (\ref{e11}) to (\ref{e12}) is an instance of the
\emph{fermionic Wick theorem}; it uses the fact that $\Psi_x$ and
$\Psi^*_y$ are linear combinations of $\psi_k$'s and $\psi^*_l$'s
respectively, together with the canonical anti-commutation
relations.

A further computation gives an explicit formula for the correlation
kernel:
$$
K(x,y)=\theta\,\frac{J_{x-\frac 12}J_{y+\frac 12}-J_{x+\frac
12}J_{y-\frac 12}}{x-y}=\sum_{k\in\Z'_{>0}}J_{x+k}J_{y+k},
$$
where $J_k=J_k(2\theta)$ are the J-Bessel functions. This is the
so-called discrete Bessel kernel that was first obtained in
\cite{Bor00b}, \cite{Joh01}.

We refer to \cite{Oko01} and \cite{Oko03} for far-reaching
generalizations of Theorem \ref{wedge}, and to \cite{Lyt02} for a
general construction of determinantal processes via representations
of the canonical anti-commutation relations corresponding to the
\emph{quasi-free states}.

\sect{Dimer models} Consider a finite planar graph $\mathcal{G}$.
Let us assume that the graph is \emph{bipartite}, i.~e. its vertices
can be colored black and white so that each edge connects vertices
of different colors. Let us fix such a coloring and denote by $B$
and $W$ the sets of black and white vertices.

A \emph{dimer covering} or a \emph{domino tiling} or a \emph{perfect
matching} of a graph is a subset of edges that covers every vertex
exactly once. Clearly, in order for the set of dimer coverings of
$\mathcal{G}$ to be nonempty, we must have $|B|=|W|$.

A \emph{Kasteleyn weighting} of $\mathcal{G}$ is a choice of sign
for each edge with the property that each face with $0\
\mathrm{mod}\ 4$ edges has an odd number of minus signs, and each
face with $2\ \mathrm{mod}\ 4$ edges has an even number of minus
signs. It is not hard to show that a Kasteleyn weighting of
$\mathcal{G}$ always exists (here it is essential that the graph is
planar), and that any two Kasteleyn weightings can be obtained one
from the other by a sequence of multiplications of all edges at a
vertex by $-1$.

A \emph{Kasteleyn matrix} of $\mathcal{G}$ is a signed adjacency
matrix of $\mathcal{G}$. More exactly, given a Kasteleyn weighting
of $\mathcal{G}$, define a $|B|\times |W|$ matrix $\K$ with rows
marked by elements of $B$ and columns marked by elements of $W$, by
setting $\K(b,w)=0$ if $b$ and $w$ are not joined by an edge, and
$\K(b,w)=\pm 1$ otherwise, where $\pm$ is chosen according to the
weighting.

It is a result of \cite{Tem61}, \cite{Kas67} that the number of
dimer coverings of $\mathcal{G}$ equals $|\det\K|$. Thus, if there
is at least one perfect matching, the matrix $\K$ is invertible.

Assume that $\det\K\ne 0$. Define a matrix $K$ with rows and columns
parameterized by the edges of $\mathcal G$ as follows:
$K(e,e')=\K^{-1}(w,b')$, where $w$ is the white vertex on the edge
$e$, and $b'$ is the black vertex on the edge $e'$. The next claim
follows from the results of \cite{Ken97}.

\begin{theo}\label{dimer} Consider the uniform measure on the dimer covers of
$\mathcal{G}$ as a random point process on the set $\X$ of edges of
$\mathcal{G}$. Then this process is determinantal, and its
correlation kernel is the matrix $K$ introduced above.
\end{theo}

The theory of random dimer covers is a deep and beautiful subject
that has been actively developing over the last 15 years. We refer
the reader to \cite{Ken08} and references therein for further
details.

\sect{Uniform spanning trees}

Let $\mathcal{G}$ be a finite connected graph. A \emph{spanning
tree} of $\mathcal{G}$ is a subset of edges of $\mathcal{G}$ that
has no loops, and such that every two vertices of $\mathcal{G}$ can
be connected within this subset.

Clearly, the set of spanning trees is nonempty and finite. Any
probability measure on this set can be viewed as a random point
process on the set $\X$ of edges of $\mathcal{G}$. We are interested
in the uniform measure on the set of spanning trees, and we denote
the corresponding process on $\X$ by $\P$.

Let us fix an orientation of all the edges of $\mathcal{G}$. For any
two edges $e=\vec{xy}$ and $f$ denote by $K(e,f)$ the expected
number of passages through $f$, counted with a sign, of a random
walk started at $x$ and stopped when it hits $y$.

The quantity $K(e,f)$ also has an interpretation in terms of
electric networks. Consider $\mathcal{G}$ with the fixed orientation
of the edges as an electric network with each edge having unit
conductance. Then $K(e,f)$ is the amount of current flowing through
the edge $f$ when a battery is connected to the endpoints $x$ and
$y$ of $e$, and the voltage is such that unit current is flowing
from $y$ to $x$. For this reason, $K$ is called the \emph{transfer
current matrix}.

This matrix also has a linear-algebraic definition. To any vertex
$v$ of $\mathcal{G}$ we associate a vector $a(v)\in \ell^2(\X)$
(recall that $\X$ is the set of edges) as follows:
$$
a(v)=\sum_{e\in\X} a_e(v) \delta_e,\qquad a_x(e)=\begin{cases}
1,&\text{if $v$ is the tail of $e$},\\
-1,&\text{if $v$ is the head of $e$},\\
0,&\text{otherwise}.
\end{cases}
$$
Then $[K(e,f)]_{e,f\in\X}$ is the matrix of the orthogonal
projection operator with image $\mathrm{Span}\{a(v)\}$, where $v$
varies over all vertices of $\mathcal{G}$, see e.~g. \cite{Ben01}.

\begin{theo}\label{USP} $\P$ is a determinantal point process  with
correlation kernel $K$.
\end{theo}

Theorem \ref{USP} was proved in \cite{Bur93}; another proof can be
found in \cite{Ben01}. The formulas for the first and second
correlation functions via $K$ go back to \cite{Kir1847} and
\cite{Bro40}, respectively. We refer to \cite{Lyo03} and references
therein for further developments of the subject.

Note that for planar graphs, the study of the uniform spanning trees
may be reduced to that of dimer models on related graphs and
\emph{vice versa}, see \cite{Bur93}, \cite{Ken00}, \cite{Ken04}.

\sect{Hermitian correlation kernels}

Let $\X$ be $\R^d$ or $\Z^d$ with the Lebesgue or the counting
measure as the reference measure. Let $K$ be a nonnegative operator
in $L^2(\X)$\footnote{An Hermitian $K$ must be nonnegative as we
want $\det[K(x_i,x_j)]\ge 0$.}. In the case $\X=\R^d$, we also
require $K$ to be locally trace class, i.~e.~for any compact
$B\subset\X$, the operator $K\cdot\mathbf{1}_B$ is trace class. Then
Lemma 2 in \cite{Sos00} shows that one can choose an integral kernel
$K(x,y)$ of $K$ so that
$$
\mathrm{Trace}(K\cdot\mathbf{1}_B)^k=\int_{B^k}K(x_1,x_2)K(x_2,x_3)\cdots
K(x_k,x_1)dx_1\cdots dx_k,\quad k=1,2,\dots
$$
For $\X=\Z^d$, $K(x,y)$ is just the matrix of $K$.

\begin{theo}\label{herm} There exists a determinantal point process on $\X$ with
the correlation kernel $K(x,y)$ if and only if $0\le K\le 1$,
i.~e.~both $K$ and $\mathbf{1}-K$ are nonnegative.
\end{theo}

Theorem \ref{herm} was proved in \cite{Sos00}; an incomplete
argument was also given in \cite{Mac75}. Remarkably, it remains the
only known characterization of a broad class of kernels that yield
determinantal point processes.

Although only Theorem \ref{bi} with $\Phi_i=\overline{\Psi}_i$,
Theorem \ref{pp}, and Theorem \ref{USP} from the previous sections
yield manifestly nonnegative kernels, determinantal processes with
such kernels are extremely important, and they are also the easiest
to analyze asymptotically, cf. \cite{Hou06}.

Let us write down the correlation kernels for the two most widely
known determinantal point processes; they both fall into the class
afforded by Theorem \ref{herm}.

The \emph{sine process} on $\R$ corresponds to the \emph{sine
kernel}
$$
K^\mathrm{sine}(x,y)=\frac{\sin\pi(x-y)}{\pi(x-y)}=\int_{-\frac
12}^\frac 12 e^{2i\pi\tau x}e^{-2i\pi\tau y}d\tau.
$$
The Fourier transform of the corresponding integral operator
$K^\mathrm{sine}$ in $L^2(\R)$ is the operator of multiplication by
an indicator function of an interval; hence $K^\mathrm{sine}$ is a
self-adjoint projection operator.

The \emph{Airy point process}\footnote{Not to be confused with the
Airy process that describes the time evolution of the top particle
of the Airy point process, see Chapter 37 of the present volume.} on
$\R$ is defined by the \emph{Airy kernel}
$$
K^\mathrm{Airy}(x,y)=\frac{Ai(x)Ai'(y)-Ai'(x)Ai(y)}{x-y}=\int_0^{+\infty}Ai(x+\tau)Ai(y+\tau)d\tau,
$$
where $Ai(x)$ stands for the classical Airy function. The integral
operator $K^\mathrm{Airy}$ can be viewed as a spectral projection
operator for the differential operator $\frac{d^2}{dx^2}-x$ that has
the shifted Airy functions $\{Ai(x+\tau)\}_{\tau\in\R}$ as the
(generalized) eigenfunctions.

\sect{Pfaffian point processes}

A random point process on $\X$ is called \emph{Pfaffian} if there
exists a $2\times 2$ matrix valued skew-symmetric kernel $K$ on $\X$
such that the correlation functions of the process have the form
$$
\rho_n(x_1,\dots,x_n)=\mathrm{Pf}\left[K(x_i,x_j)\right]_{i,j=1}^n,\qquad
x_1,\dots,x_n\in\X,\quad n=1,2,\dots
$$
The notation $\mathrm{Pf}$ in the right-hand side stands for the
Pfaffian, and we refer to \cite{deB55} for a concise introduction to
Pffafians.

Pfaffian processes are significantly harder to study than
determinantal ones. Let us list some Pffafian analogs of the
statements from the previous sections.

\begin{itemize}

\item A Pfaffian analog of the Fredholm determinant formula
(\ref{gen}) for the generating functional can be found in Section 8
of \cite{Rai00}.

\item A Pfaffian analog of the Eynard-Mehta theorem is
available in \cite{Bor05}.

\item Pfaffians can be used to enumerate nonintersecting paths
with free endpoints, see \cite{Ste90}. This leads to combinatorial
examples for the Pfaffian Eynard-Mehta theorem.

\item Pfaffian L-ensembles and conditional L-ensembles are treated
in \cite{Bor05}.

\item Fermionic Fock space computations leading to a Pfaffian point process
were performed in \cite{Fer04} and \cite{Vul07}.

\item Pfaffians arise in the enumeration of dimer covers of planar
graphs that are not necessarily bipartite, see \cite{Tem61},
\cite{Kas67}.

\end{itemize}

\ \\
{\sc Acknowledgements}: This work was supported by NSF grant
DMS-0707163.

\end{document}